\documentclass[11pt]{article}

\usepackage{amssymb, amsbsy}
\usepackage{amsmath, amsthm}
\usepackage{amsfonts}

\usepackage{graphicx}
\usepackage{epsfig}

\newcommand{\aac}{\`a }

\newcommand{\FF}{{\mathbb F}}

\newcommand{\A}{{\alpha}}

\newcommand{\LL}{{\ldots}}
\newcommand{\E}{{\varepsilon}}

\newcommand{\OO}{{\Omega}}

\newcommand{\HH}{H_{\alpha}}
\newcommand{\FFqn}{\mathbb{F}_{q^n}}
\newcommand{\FFqnr}{\mathbb{F}_{q^{nr}}}
\newcommand{\Sqn}{S_{q^{n}}}
\newcommand{\SSS}{\mathbb{S}}

\long\def\symbolfootnote[#1]#2{\begingroup%
\def\thefootnote{\fnsymbol{footnote}}\footnote[#1]{#2}\endgroup}
\newtheorem{thm}{Theorem}[section]
\newtheorem{theorem}[thm]{Theorem}
\newtheorem{corollary}[thm]{Corollary}

\newtheorem{claim}[thm]{Claim}

\newtheorem{fact}[thm]{Fact}

\newtheorem{proposition}[thm]{Proposition}

\newtheorem{definition}[thm]{Definition} \theoremstyle{definition}
\theoremstyle{remark}
\newtheorem{remark}[thm]{Remark}

\begin{document}

\title{Permutation equivalent maximal irreducible Goppa codes}

\author{  Francesca  Dalla Volta\footnote{francesca.dallavolta@unimib.it, Dipartimento di Matematica e applicazioni 
Universit\aac di Milano Bicocca },  Marta Giorgetti\footnote{ marta.giorgetti@uninsubria.it,  Dipartimento di Fisica e Matematica
Universit\aac dell'Insubria, Como}, Massimiliano Sala\footnote{ msala@bcri.ucc.ie,  Dipartimento di Matematica
  Universit\aac di Trento}}
\bibliographystyle{amsalpha}

\maketitle

\begin{abstract}
\noindent We consider the problem of finding the number of permutation non-equivalent classical irreducible maximal Goppa codes
having fixed parameters $q,\,n$ and $r$ from a group theory point of view.


\end{abstract}

\section{Introduction}

The study of classical Goppa codes is
important since they are a very large class of codes, near to random
codes. They are easy to ge\-nerate and possess an interesting
algebraic structure. For these reasons they are used in McEliece's public key
cryptosystem  \cite{CGC-cd-art-McEli1}. This cryptosystem is based on the difficulty to find a generator matrix of a Goppa code when a "scrambled" of it is known.

In this paper we consider the problem of finding an upper
bound for the number of permutation non-equivalent irreducible maximal Goppa codes. This question was considered by several authors (see for example \cite{CGC-cd-art-chen}, \cite{berger1}, \cite{berger2}, \cite{CGC-cd-art-patjohn1}, \cite{CGC-cd-phdthesis-marta}).
In Section 3 we briefly recall these approaches. In particular, we describe the action of a group $FG$ isomorphic to $A\Gamma L(1,q^n)$ on the $q^n$ columns of a suitable parity-check matrix $H_{\alpha}$. This induces on maximal irreducible Goppa codes the same action which arises from \cite{CGC-cd-phdthesis-john}. This action does not describe exactly the orbits of Goppa codes,
    since in some cases the number of permutation non-equivalent
    Goppa codes is less than the number of orbits of $FG$.
The group $FG$ acts faithfully on columns of $H_{\A}$, so that it may be seen as a subgroup of the symmetric group $\Sqn $. It seems interesting to study if there is a proper subgroup of $\Sqn $ containing $FG$, acting on the set $\Omega $ of classical irreducible maximal Goppa codes of fixed parameters, and giving on $\OO$  exactly the orbits of permutation equivalent codes. In order to consider this problem, we analyze the subgroups of $\Sqn $ containing $FG$ and
in Section 4 we find that there exists exactly one maximal subgroup $M$, isomorphic to $AGL(nm,p)$ of $\Sqn $ ($A_{q^n}$) containing $FG$ ($q=p^m)$. This suggests that one could  consider the action of this $M$ on codes to reach the right bound.

We are grateful to Andrea Caranti, Andrea Lucchini, John A. Ryan and Patrick Fitzpatrick for helpful discussions on this subject.

\section{Preliminaries}
In this section we fix some notation and we recall some basic concepts about linear codes and in particular about Goppa codes. Our main references are \cite{MR1996953} for coding theory and \cite{CGC-alg-book-cameron1} for  group theory.

We denote by $\FF_q$ the finite field with $q$ elements, where $q=p^m$ is a power of a prime $p$; let $N$, $k$, $n$ and $r$ be natural numbers, $k \leq N$.
We consider two extensions of
$\FF_q$, of degree $n$ and $nr$, $\FF_{q^{n}}$ and
$\FF_{q^{{n}r}}$ respectively; $\FF_{q^{n}}[x]$ denotes
the polynomial ring over $\FF_{q^{n}}$
and $\E$ is a primitive element of $\FF_{q^{n}}$, $\FF_{q^{n}}^*=\langle\E\rangle$. We refer to the vector space of dimension $N$ over $\FF_q$ as to $(\FF_q)^N$.

In the following
if $H$ is an $(N-k)\times N$ matrix with entries in $\FF_{q}$
and rank equal to $N-k$, the set $C$ of all
vectors $c \in(\FF_q)^{N}$ such that $H c^T=0$ is an $(N,k)$ {\it
linear} code over $\FF_q$, of {\it length}  $N$ and
{\it dimension} $k$, i.e. a subspace of $(\FF_q)^N$ of dimension $k$. The elements of $C$ are called  {\it codewords} and
matrix $H$ is a {\it parity-check matrix} of $C$.  Any $k \times N$ matrix $G$ whose rows form a
vector basis of $C$ is called a {\it generator matrix} of  $C$.

\begin{definition}
Let $E/K$ be a field extension. A linear code $C$ is called a \textbf{ subfield subcode} if $C$ is obtained as the restriction
to $K^n$ of a linear subspace $L$ of $E^n$.
\end{definition}

By abuse of notation  we call parity-check matrix also a matrix $H$ with entries in $E$  such that $Hc^T=0$ for all $c\in C$. According to this assumption, $H_1$ and $H_2$ may be parity-check matrices for the same code even if their entries are in different extension fields or they have different ranks.

\begin{definition}[\cite{MR1996953}]\label{def:equivalenze:forte}
Let $C_1$ and $C_2$ be two linear codes  over $\FF_q$ of length $N$, let $G_1$ be a generator matrix of $C_1$. Codes $C_1$ and $C_2$ are
\textbf{permutation equivalent} provided there is a permutation $\sigma \in S_N$
of coordinates which sends $C_1$ in $C_2$. Thus $C_1$ and $C_2$ are permutation equivalent provided there is a permutation matrix $P$ such that $G_1P$ is a generator matrix for $C_2$.
They are \textbf{monomially equivalent} provided there is a monomial matrix $M$ so that $G_1M$ is a generator matrix for $C_2$ and
\textbf{equivalent} provided there is a monomial matrix $M$ and an automorphism $\gamma$ of the field $\FF_q$ so that $C_2=C_1M\gamma$.

\end{definition}

If code $C_2$ is
permutation equivalent to $C_1$ with parity-check matrix $H_1$, we can obtain a parity-check
matrix $H_2$  for $C_2$ by permuting columns of $H_1$ (and viceversa).

%
%
%

\begin{definition}\label{goppa:def}
Let $g(x)=\sum g_i x^i\in \mathbb{F}_{q^{n}}[x]$ and let
$L=\{\E_1,\E_2,\LL,\E_N\}$ denote a subset of elements of
$\mathbb{F}_{q^{n}}$ which are not roots of $g(x)$. Then the
\textbf{Goppa code} $\mathcal{G}(L,g)$ is defined as the set of all
vectors $c=(c_1,c_2,\LL,c_N)$ with components in $\mathbb{F}_q$
which satisfy the condition:
\begin{equation}\label{goppa:eq}
\sum_{i=0}^{N}\frac{c_i}{x-\E_i}\equiv 0 \mod\,g(x).
\end{equation}
\end{definition}
Usually, but now always, the set $L=\{\E_1,\E_2,\LL,\E_N\}$ is taken
to be the set of all elements in $\mathbb{F}_{q^{n}}$ which are
not roots of the Goppa polynomial $g(x)$.
In this case the Goppa
code is said \textit{maximal}.
If the degree of $g(x)$ is $r$, then the Goppa code is called a Goppa
code of degree~$r$. It is easy to see
(\cite{CGC-cd-book-macwilliamsI}) that a parity-check matrix for
$\mathcal{G}(L,g)$ is given~by
$$
H=\left(
        \begin{array}{cccc}
            \frac{1}{g(\varepsilon_1)} & \frac{1}{g(\varepsilon_2)} & \LL & \frac{1}{g(\varepsilon_N)} \\
            \frac{\varepsilon_1}{g(\varepsilon_1)} & \frac{\varepsilon_2}{g(\varepsilon_2)} & \LL & \frac{\varepsilon_N}{g(\varepsilon_N)} \\
            \vdots & \vdots & \vdots & \vdots \\
\frac{\varepsilon_1^{r-1}}{g(\varepsilon_1)}
&\frac{\varepsilon_2^{r-1}}{g(\varepsilon_2)} & \LL &
\frac{\varepsilon_N^{r-1}}{g(\varepsilon_N)}\\
\end{array}
    \right).
$$
Note that the code $C=\ker H$ is a subspace of $(\FF_{q^n})^N$ and the Goppa code $\mathcal{G}(L,g)$ is its subfield subcode on $\FF_q$.


\begin{definition}
A Goppa code $\mathcal{G}(L,g)$ is called \textbf{irreducible} if $g(x)$
is irreducible over $\mathbb{F}_{q^{n}}$.
\end{definition}
In the following by Goppa code we mean maximal irreducible classical Goppa code of degree $r$,  so that $N=q^n$. By Definition \ref{goppa:def}, a
vector
$c=(c_1,c_2,\LL,c_{q^{n}})\in(\mathbb{F}_q)^{q^{n}}$ is a
codeword of $\mathcal{G}(L,g)$ if and only if it satisfies
(\ref{goppa:eq}).
If $\alpha$ is any root of $g(x)$, $\A\in \FFqnr$, then
$g(x)=\prod_{i=0}^{r-1}(x-\alpha^{q^{{n}i}})$ and (\ref{goppa:eq}) is equivalent to the $r$ equations
\begin{equation}\label{seconda}
\sum_{i=1}^{q^{n}}\frac{c_i}{\A^{q^{{n}
j}}-\E_i}=0,\,\quad 0\leq j\leq r-1.
\end{equation}

Hence $\mathcal{G}(L,g)$ is completely described by any root $\alpha$ of
$g(x)$ and we may denote this code by $\mathcal{C}(\alpha)$. From (\ref{seconda}) we easily  get a
parity-check matrix $H_{\A}\in Mat_{1\times q^n}(\FFqnr)$ for $\mathcal{C}(\alpha)$ (see \cite{CGC-cd-art-chen}):
\begin{equation}\label{hriga}
H_{\A}=\left(
\begin{array}{cccc}
\frac{1}{\alpha-\E_1}, & \frac{1}{\alpha-\E_2},  & \LL, &
\frac{1}{\alpha-\E_{q^{n}}} \\
\end{array}%
\right).
\end{equation}
It is important
to stress that by using parity-check matrix $H_{\A}$ to define
$\mathcal{C}(\alpha)$ we  implicitly fix an order in $L$. So, we set $$L=\{\E,\E^2,\LL,\E^{q^{{n}-1}},\E^{-\infty}\},
$$
where $ \E^{-\infty}=0$, $\E_i=\E^i$  and  the matrix $H_{\A}$ is
$$
H_{\A}=\left(
\begin{array}{ccccc}
\frac{1}{\alpha-\E}, & \frac{1}{\alpha-\E^2},  & \LL, &
\frac{1}{\alpha-1}, &\frac{1}{\alpha}\\
\end{array}%
\right).
$$
We observe that the Goppa code $C(\A)$ is the subfield subcode of codes having as parity-check matrices both  $H$ and $H_{\A}$.
Moreover, there exist matrices having  structure different from $H$ and $H_{\A}$, which are parity-check matrices for $C$.

We denote by $\OO=\OO(q,n,r)$  the set of Goppa codes, with fixed
parameters $q,\,n,\,r$.

In the following an action on set $\SSS$ is considered, where $\mathbb{S}=\mathbb{S}(n,r)$ is composed of all elements in
$\mathbb{F}_{q^{{n} r}}$ of degree $r$ over
$\mathbb{F}_{q^{n}}$.

\section{Three actions on $\OO$}
In this section we briefly present semiaffine actions introduced in \cite{berger3} and in \cite{CGC-cd-art-patjohn2}. These actions have degrees $\frac{|\SSS|}{r}$ and $|\SSS|$ respectively. Moreover we consider an action of the group $A\Gamma L(1,q^n)$
on entries of parity-check matrix of type $\HH$. This time the degree is $q^n$.

\vspace{0,5cm}

In \cite{berger3}, the author works directly on polynomials by studying automorphism groups of several classes of codes.

If $\psi\in AG L (1,{q^n})$, $\psi (z)= az+b,\,\, a,b\in \FFqn,\, a\neq 0 $, he defines
$$
g^{\psi}(x)=\sum_{i=0}^rg_i(ax+b)^i.
$$
The map $\psi$ acts also on set $L=\FFqn$, $\FFqn^*=\langle\varepsilon\rangle$, by
$$
L^{\psi}=\left(\E^{\psi^{-1}},\LL, (\E^{q^n-1})^{\psi^{-1}}, (\E^{-\infty})^{\psi^{-1}}\right).
$$
The code $\mathcal{G}(L,g)^{\psi}=\mathcal{G}(L^{\psi},g^{\psi})$ is said the \textit{conjugate} of code $\mathcal{G}(L,g)$ by $\psi$.

\begin{proposition}{\cite{berger3}}\label{prop_berger}
The Goppa codes are invariant by conjugation under the affine group $AGL(1,{q^n})$, i.e. $\mathcal{G}(L,g)^{\psi}=\mathcal{G}(L,g)$ for all $\psi$ such that $\psi(z)=az+b, a,b \in \FF_{q^n}, a\neq0$.
\end{proposition}

We get
\begin{corollary}
Goppa codes $\mathcal{G}(L,g^{\psi})$ is equivalent to Goppa code $\mathcal{G}(L,g)$.
\begin{proof}
We known that $\mathcal{G}(L,g)^{\psi}=\mathcal{G}(L^{\psi},g^{\psi})$ and from Proposition \ref{prop_berger} $\mathcal{G}(L,g)^{\psi}=\mathcal{G}(L,g)$.
From Definition \ref{def:equivalenze:forte} it follows that $\mathcal{G}(L,g^{\psi})$ is equivalent to $\mathcal{G}(L^{\psi},g^{\psi})=\mathcal{G}(L,g)$ so $\mathcal{G}(L,g)$ is equivalent to $\mathcal{G}(L,g^{\psi})$.
\end{proof}
\end{corollary}

More generally, if $\psi \in A\Gamma L(1, q^n)$, $\psi (z)=az^{q^t}+b$, with $a,b \in \FFqn,\, a\neq 0$ ant $t\in \{0,\LL,n-1\}$, we define
\begin{equation}\label{psi}
g^{\psi}(x)=\sum_{i=0}^rg_i(ax^{q^t}+b)^i\;\;
\end{equation}
Equation (\ref{psi}) suggests to consider an action $\sigma$ on $\mathbb{P}\subseteq \FFqn[x]$, where $\mathbb{P}$ is the set of irreducible polynomials of degree $r$. For $g\in \mathbb{P}$, $g^{\sigma(\psi)}$ is the unique  polynomial $f$ of degree $r$ such that $g(\A)=0$ if and only if $f(\beta)=0$ for  $\beta=\left(\frac{\A-b}{a}\right)^{q^{nr-t}}$ (note  $g^{\sigma(\psi)}\in \mathbb{P}$). 

Indeed, if $g(x)=\sum_{i=0}^r g_i x^i$, there exist $\bar{g}_i, \forall i=1,\ldots,r,\;  \bar{a},\; \bar{b}$ such that $\bar{g}_i^{q^t}=g_i, \forall i=1,\ldots,r,\; \bar{a}^{q^t}=a,\; \bar{b}^{q^t}=b$ so that
$$\sum_{i=0}^r \bar{g}_i^{q^t} (\bar{a}^{q^t}x^{q^t}+\bar{b}^{q^t})^i=\left(\sum_{i=0}^r \bar{g}_i (\bar{a}x+\bar{b})^{i}\right)^{q^t}=f(x)^{q^t}.$$
    It is immediate to recognize that $g(\A)=0$, for $\A\in\SSS$, if and only if $f(\beta)=0$, with $\beta=\left(\frac{\A-b}{a}\right)^{q^{nr-t}}\in \SSS$.

With similar arguments used for  Proposition \ref{prop_berger}, one gets

\begin{proposition}\label{ber:1}
The Goppa codes are invariant by conjugation under the semiaffine group $A\Gamma L(1,{q^n})$, i.e. $\mathcal{G}(L,g)=\mathcal{G}(L^{\psi}, f)$, where $f=g^{\sigma(\psi)}\in \mathbb{P}$.
\end{proposition}

\begin{corollary}
Goppa codes $\mathcal{G}(L,g)$ is equivalent to Goppa code $\mathcal{G}(L,f)$.
\end{corollary}
\vspace{0,5cm}
In \cite{CGC-cd-phdthesis-john}  the same action on $\OO$ is obtained considering an action on $\SSS$ of an "affine" group $T=AGL(1,q^n)\langle\sigma\rangle$, where $\sigma$ is defined as $\sigma: x\rightarrow x^{q}$; the group $\langle \sigma \rangle$ has order $nr$.
The main result is the following:

\begin{theorem}\cite{CGC-cd-phdthesis-john}\label{alpha-beta}
If $\A,\,\beta\in \mathbb{S}$ are related as it follows
\begin{equation}\label{relazionealphabeta}
\beta=\zeta \A^{q^i} + \xi
\end{equation}
for some $\zeta, \xi \in \FF_{q^{n}},\, \zeta  \neq 0, i= 1\LL nr$, then
$C(\A)$ is equivalent to $C(\beta)$.
\end{theorem}

 Orbits over $\SSS$ give orbits on $\OO$.

\begin{fact}
The above actions on $\SSS$ and on $\mathbb{P}$ create the same orbits on $\OO$.
\end{fact}
\begin{proof}
Let $\A\in\SSS$ be a root of $g(x)\in\mathbb{P}$. Let $\beta=\zeta\A^q+\xi\in \SSS$. There exists an irreducible polynomial $g_1\in \FF_{q^n}[x]$,  such that $g_1(\beta)=0$. From Proposition \ref{ber:1} we get that the orbit  $\A^{FG}=\{ t(\beta ),\,  t \in T\}$) induces on $\OO $ the same  orbit than $g^{T}=\{t( g), t\in T\}$).
\end{proof}

The work in \cite{berger3} is mainly directed to the study of automorphism group of a given code;
\cite{CGC-cd-phdthesis-john} is deeply interested in counting the number of non-equivalent Goppa codes.

In \cite{CGC-cd-phdthesis-john} the exact number of orbits on $\SSS$ is given. Unfortunately, several examples are exhibited where the number of orbits $T$ on $\SSS$ is greater than the number of non-equivalent Goppa codes.

We introduce an action on columns of $H_{\A}= \left(\frac{1}{\alpha-\E},  \frac{1}{\alpha-\E^2},   \LL,\frac{1}{\alpha-1}, \frac{1}{\alpha}\right)$ , which induces the same orbits on $\OO$ than $T$.
We state the results and give a sketch of the proofs.
For more details see \cite{CGC-cd-phdthesis-marta}.

Let us consider the subgroup $FG\simeq A\Gamma L(1, q^n)\leq S_{q^n}$ in its natural action on points of $\FF_{q^n}$. If $\psi \in FG$, then $\psi(x)=ax^{q^i}+b$, where $a,b \in \FFqn$, $a\neq 0$ and $i=1,\LL,n$.
Since each entry (column) of $\HH$ is uniquely determined by an element of $\FFqn$, $A\Gamma L(1, q^n)$ realizes a permutation of $\HH$ entries given by:
$$
\left(\frac{1}{\A-\varepsilon}\right)^{\psi}=\frac{1}{\A-\varepsilon^{\psi}}.
$$
Writing $FG$ we mean  $F=AGL(1,q^n)$ and $G$  the automorphism group of $\FFqn$ over $\FF_q$.

The matrices $H_{\beta}$ and $\zeta H_{\beta}^{q^j}$ are parity-check matrices for the same Goppa code $C(\beta)$.
We characterize the permutations mapping $H_{\A}$ into $\zeta H_{\beta}^{q^j}$ in the following proposition.

\begin{proposition}
Let $H_{\A}$ and $H_{\beta}$ be parity-check matrices for Goppa codes $C(\A)$ and $C(\beta)$. If there exists a
permutation $\rho\in S_{q^{{n}}}$, such that
$\rho(H_{\A})=\zeta (H_{\beta})^{q^j}$, for some $\zeta, \beta \in
\FF_{q^{{n}}}, \zeta \neq 0$ and $j=1,\LL,{{n}}r$, then $\rho
\in {F}{G}$.
\begin{proof}
We consider 
$$
\rho (H_{\A})=H'_{\A}=\left( \frac{1}{\A-\E^{i_1}},
\frac{1}{\A-\E^{i_2}}, \cdots ,\frac{1}{\A-\E^{i_{q^{{{n}}}-1}}},
\frac{1}{\A-\E^{i_{-\infty}}} \right),
$$
where $i_j=\rho(j)$ and
matrix $\zeta (H_{\beta})^{q^j}$:
$$
\zeta (H_{\beta})^{q^j}  = \left( \frac{\zeta}{\beta^{q^j}-\E^{q^j}},  \frac{\zeta}{\beta^{q^j}-\E^{2q^j}},  \ldots ,  \frac{\zeta}{\beta^{q^j}-1}, \frac{\zeta}{\beta^{q^j}}\right).
$$
Suppose $\zeta=1$ and ${j=1}$ so that $ \forall
t \in \{1,2,\LL,q^n\}$ we have $
\frac{1}{\beta^{q}-\E^{tq}}=\frac{1}{\A-\E^{i_t}}$ and then $\A-
\beta^q=\E^{i_t}-\E^{tq}$.

If $\A-\beta^q=0$,  $\rho$ is the permutation induced by the Frobenius map $\sigma $, since
$\E^{i_t}=\E^{tq}$; it follows that
$$ \rho(t)= \left\{
\begin{array}{ll}
tq & \mbox{if } t=1,2,\ldots,q^{{n}}-1\\
-\infty & \mbox{if } t=-\infty
\end{array}
\right.
$$
and $\rho={\sigma}.$

If $\A-\beta^q \neq 0$, as above
$\A-\beta^q\in\mathbb{F}_{q^{{n}}}$ so that $\A-\beta^q=\E^k$
for some $k\in \{1,\LL,q^{{n}}-1\}$ and then permutation $\rho
\in {F}{G}$;  explicitly it acts as:
$$
\rho(t)= \left\{
\begin{array}{ll}
i_t= tq +f_k(\E) & \mbox{if } t=1,2,\ldots,q^{{n}}-1\\
i_t= k & \mbox{if } t=-\infty
\end{array}
\right.
$$
where $i_t$ is such that $\E^{i_t}=\E^{tq}+\E^k$, and $f_k(\E)$ is a function depending on representation of $  \FFqn$.

If $\zeta \in \mathbb{F}_{q^{{n}}}^*, {\zeta \neq 1}$
and ${j=1}$, then $\zeta=\E^l$ for some $l \in
\{1,\cdots,q^{{n}}-2\}$. With same arguments used in the
previous step,  we get
$$
\zeta \A -\zeta \E^{i_t}= \beta^q- \E^{tq} \implies
\zeta\E^{i_t}=\zeta\A - \beta^q+ \E^{tq} \implies   \E^{i_t}=\A -
\zeta^{-1}\beta^q+ \E^{tq-l}.
$$
Again $\A - \zeta^{-1}\beta^q\in \FF_{q^{n}}$;  then there is $h
\in \{1,\LL,q^n\}$ so that $\E^{i_t}=\E^{h}+\E^{tq-l}$, and
$$\rho(t)= \left\{
\begin{array}{ll}
i_t= tq-l +f_h(\E) & \mbox{if } t=1,2,\ldots,q^{{n}}-1\\
i_t=h  & \mbox{if } t=-\infty
\end{array}
\right.
$$
where $\E^{i_t}=\E^{tq-l} +\E^h$ and $f_h(\E)$ depending on the representation of $\FF_{q^{n}}$. Concluding $\rho={\tau_{k}
\mu_{\zeta^{-1}} \sigma}$; here $\tau_k$ is the translation defined by $\tau_k:x\rightarrow x+\E^k$,  $\mu_{\zeta }$ is the map   $\mu_{\zeta}:x \rightarrow \zeta x$ and $\sigma $ is the Frobenius map; this proves  $\rho \in FG$.

Finally, if ${j \neq 1}$ we have:
$
\frac{\zeta}{\beta^{q^j}-\E^{tq^j}}=\frac{1}{\A-\E^{i_t}}
$
and $ \zeta \A -\zeta \E^{i_t}= \beta^{q^j}- \E^{tq^j}$. As
$\zeta=\E^l$ for some $l\in\{1,\LL,q^n-2\}$,  we gain:
$$
\E^{i_t}=\A-\E^{-l}\beta^{q^j}+ \E^{tq^j-l}.
$$
So there is $v \in \{ 1,\LL,q^n\}$ such that $ \A -
\E^l\beta^{q^j}=\E^v$ and $\E^{i_t}=\E^{v}+\E^{tq^j-l}.$
Permutation $\rho$ is:
$$
\rho(t)= \left\{
\begin{array}{ll}
i_t= tq^j-l +f_{v}(\E) & \mbox{if } t=1,2,\ldots,q^{{n}}-1\\
i_t=v  & \mbox{if } t=-\infty
\end{array}
\right.
$$
where $\E^{i_t}=\E^{tq^j-l}+\E^v$ and $f_v(\E)$ depends on the representation of $\FF_{q^n}$. Concluding $\rho={\tau_{v}
\mu_{\zeta^{-1}} \sigma^j}.$ Clearly in all cases $\rho\in
{F}{G}$.
\end{proof}
\end{proposition}

\begin{corollary}
Let $H_{\A}$ and $H_{\beta}$ be parity-check matrices for Goppa codes $C(\A)$ and $C(\beta)$. If there exists a
permutation $\rho\in S_{q^{{n}}}$, such that $\rho(H_{\A})=
\zeta H_{\beta}$, for same $\zeta, \beta \in \FF_{q^{n}}, \zeta
\neq 0$, then $\rho\in {F}$.

\end{corollary}
\section{Maximal subgroups}

\def\Aqn{A_{q^n}}
The action of $A\Gamma L(1,q^n)$ does not reach the exact number of non-equivalent maximal Goppa codes. So we look for  maximal subgroups of $S_{q^n}$ containing a fixed $A\Gamma L(1,q^n)=FG$.

\begin{theorem}[\cite{CGC-alg-book-cameron1}]\label{th_cameron1}
A maximal subgroup of $S_{q^n}$ is one of the following:
\begin{enumerate}
\item intransitive, $S_k \times S_l$, $k+l=q^n$;
\item transitive imprimitive: the wreath product $S_k \,Wr\, S_l$ in the  standard action,
$kl=q^n$;
\item primitive non-basic, the wreath product $S_k \,Wr\, S_l$ in the product action,
$k^l=q^n, k\neq 2$;
\item affine $AGL(d,p), p^d=q^n$;
\item diagonal, $T^k.(Out(T)\times S_k$), $T$ non abelian
simple, $|T|^{k-1}=q^n$; here $Out(T)$ denotes, as usual, the factor group $\frac{Aut(T)}{T}$.
\item almost simple, that is an automorphism group $G$ of a finite non abelian simple group $S$, $S\leq G\leq Aut (S)$.
\end{enumerate}
A maximal subgroup of the alternating group is the intersection of one of these groups with the alternating group.

\end{theorem}
\begin{remark}
We explicitly observe that for $p$ even, $d\geq 3$, the group $AGL(d,p)$ is actually contained in the alternating group $A_{p^d}$.
It is sufficient to realize that,  in this case, the translations are product of $2^{d-1}$  cycles of length $2$, as well as the transvections are product of $2^{d-2}$  cycles of length $2$.  As the transvections generate the general linear group $GL(p,2)$, $AGL(p,2)$ is contained in $A_{2^d}$.
\end{remark}

\begin{proposition}\label{AGL}
 $FG $ is contained in  $\Aqn $ if and only if  $q$ is even.
\end{proposition}
\begin{proof}
The thesis follows from the following result.
\begin{claim}\label{Lee}\cite{CGC-alg-Li}
Let $X$ be a primitive permutation group of degree $n$. Then $X$ contains an abelian regular subgroup $G$ if and only if either
\begin{enumerate}
\item [i)] $X \leq AGL(d,p)$, where $p$ is a prime, $d \geq 1$ and $n = p^d$; or
\item [ii)] $X=(\tilde{T_1} \times \LL \times \tilde{T_l})\cdot O \cdot
P$, $G={G_1} \times \LL \times {G_l}$ where $n={m}^l$, $l\geq1$, $G_i
<\tilde{T}_i$, with $|G_i|={m}$, $\tilde{T}_1\cong\LL\cong
\tilde{T}_l$, $O\leq
\mathrm{Out}(\tilde{T}_1)\times\LL\times\mathrm{Out}(\tilde{T}_l)$,
$ P$ is a transitive permutation group of degree $l$ and one of the following holds:
\begin{itemize}
\item[(a)] $(\tilde{T}_i,G_i)=(PSL(2,11),\mathbb{Z}_{11})$,
$(M_{11},\mathbb{Z}_{11})$, $(M_{12},\mathbb{Z}_{2}^2\times
\mathbb{Z}_3)$, $(M_{23},\mathbb{Z}_{23})$ ($M_i$ are the Mathieu groups);
\item[(b)] $\tilde{T}_i=PGL(d,q)$ e $G_i=\mathbb{Z}_{\frac{q^d-1}{q-1}}$
is a Singer group;
\item[(c)] $\tilde{T}_i=P\Gamma L(2,8)$ and
$G_i=\mathbb{Z}_9\nleq PSL(2,8)$;
\item[(d)] $\tilde{T}_i=S_{m}$ or $A_{m}$ and $G_i$ is an abelian group of order ${m}$.
\end{itemize}
\end{enumerate}
  \end{claim}
  Take $X=FG $. $FG$ contains the subgroup $A$ of translations, $A=\{\tau_{\E}:x\rightarrow x+\E\}$, so that $FG$ is contained in $N_{S_{q^n}}(A)=AGL(nm,p)$. By the above remark, if $p=2$, the group  $FG $ is contained in $ \Aqn$. If $p$ is  odd,  then the element $\mu_{\E}:x\rightarrow \E x$ belongs to $FG$ and it is odd, as its order is $q^n-1$ (recall that an element of order $q^n-1$ is said a Singer cycle); this proves that   $FG $ (and $AGL(nm,p)$) is not a subgroup of $\Aqn $.
\end{proof}

\begin{thm} Let $G=\Aqn $  if $q=2^m$,  $G= \Sqn $ for $q$ odd.
If $M$ is a maximal subgroup of $G $ containing $ F G$, then  $M$ is isomorphic to the affine group $AGL(nm,p)$ . Moreover, there is exactly one maximal subgroup containing $  F G$.

\end{thm}

\begin{proof}
As $FG $ is a primitive  $2$-transitive group of $G$,  a  maximal subgroup $M$ of $\Sqn $ containing  $ F G$,  is an almost simple group or it is isomorphic to the affine group $AGL(nm,p)$ (\cite{CGC-alg-book-cameron1}).
In the proof of Proposition \ref{AGL} we have seen that  $FG $ is contained in $AGL(mn,p)$. We prove that it is not contained in an almost simple group. By contradiction, let $M$ be an automorphism group of a simple non abelian group $S$,  $S\leq M\leq Aut (S)$. If $M$ contains $FG$, the stabilizer of a point $\omega$ in  $\FF _{q^n}$ has index $q^n=p^{nm}$. As $S$ is normal in $M$, $S$ is transitive on $\FF _{q^n}$, so that we are reduced to consider subgroups of prime power index in $S$. These were  described by Guralnick and
 for the reader's sake  we write  the main result of \cite{CGC-alg-art-Guralnick1}.
\begin{claim}[\cite{CGC-alg-art-Guralnick1}]\label{Gural}
Let $G$ be a nonabelian simple group with $H<G$ and $[G:H]=p^d=q^n$, $p$ prime. One of the following holds.
\begin{enumerate}
\item $G=A_{q^n}$ and $H\cong A_{q^n-1}$;
\item $G=PSL(s,t)$ and $H$ is the stabilizer of a line or hyperplane. Then $[G:H]=\frac{t^s-1}{t-1}=q^n$ (Note $s$ must be prime);
\item $G=PSL(2,11)$ and $H\cong A_5$;
\item $G=M_{23}$ and $H\cong M_{22}$ or $G=M_{11}$ and $H\cong
M_{10}$;
\item $G=PSU(4,2)\cong PSp(4,3)$, $H$ is the parabolic subgroup of
index 27.
\end{enumerate}
\end{claim}

 Cases 3, 4, 5,  are easily ruled out, as $p^{mn}$ is neither a prime number, nor $27$. Similarly,  case 1 is ruled out when  $p$ is  odd, as, in this case, the element $\mu_{\E}$ is odd.  If $p=2$, then  $FG $ is actually contained  in $M\simeq \Aqn $.
So, we are left with Case 2. Here, we use Claim \ref{Lee}. $X$ satisfies condition ii), with
$$ X=S=PSL(s, t), \qquad [S:H]=\frac{s^t-1}{t-1}=p^{nm},\, \, l=1.$$
  and it is easy to see that it is not the case.

\noindent Now, we prove that there is exactly one subgroup isomorphic to $AGL(nm,p)$ containing $FG $.

Let $q$ be odd:
in $\Sqn $  there is exactly one  conjugacy class of maximal subgroups of this type (see for example \cite{conjugacyaffine}).
So, let $FG \leq M\simeq AGL(nm,p)$, where the normal subgroup of the translation of $M$ is exactly the translation
group $A$ of $FG$ (\cite{CGC-alg-book-hupper1}). The element $\mu_{\E} $ generates a Singer subgroup; it is well known that the Singer cycles are conjugated in $M$; from the knowledge of the overgroups of a Singer cycle \cite{kantor_singer},
\cite{CGC-cd-art-britnellguralnickmaroti}, one easily proves that also  the normalizers of  Singer cycles contained in $M$ are  conjugate in $M$. It follows that if  $FG ^g$, $g\in \Sqn$ is contained in $M$, there exists an element $m\in M$, such that $FG ^g=FG ^m$. So, if $s$ denotes the number of the subgroups of $M$ containing $FG$, we get:
$$[\Sqn :N_{\Sqn}(FG)] = \frac{[\Sqn :M]\cdot [M:N_{M}(FG)]}{s};$$
now, from \cite{CGC-alg-Li} one gets $N_{\Sqn}(FG)\leq M$, so that $s=1$.

Now, suppose $q$ is even. $FG\leq \Aqn$ and in $\Aqn $ the conjugacy class of $S_{q^n}$ subgroups   which are isomorphic to $AGL(nm,p)$. In $A_{q^n}$ $AGL(nm,p)$ splits into two classes so that also the class of Singer cycles splits into two different classes. Same argument used for the odd case leads to the result.
\end{proof}


\begin{thebibliography}{4321}
\bibitem{berger3} {Berger, Thierry P.}, \textit{On the cyclicity of {G}oppa codes, parity-check subcodes of
 {G}oppa codes, and extended {G}oppa codes}, {Finite Fields and their Applications}, {6}, {2000}.


\bibitem{berger1} {Berger, Thierry P.}, \textit{{Cyclic alternant codes induced by an automorphism of a {GRS}
 code}}, {Finite fields: theory, applications, and algorithms (Waterloo, ON, 1997)}.

\bibitem{berger2} {Berger, Thierry P. and Charpin, P.}, \textit{{The permutation group of affine-invariant extended cyclic
 codes}}, {IEEE Trans. Inform. Theory}, {42}, {1996}.



\bibitem{CGC-cd-art-britnellguralnickmaroti}
{Britnell, John R. and  Evseev, Anton and Guralnick, Robert M. and  Holmes, Petra E. and Maroti, Attila},
    \textit{{Sets of elements that pairwise generate a linear group}},
  {preprint {\tt www-circa.mcs.st-and.ac.uk}},
{2006}.

\bibitem{CGC-alg-book-cameron1}  {Cameron, P. J.}, \textit{{Permutation groups}},  {Cambridge University Press},  {1999},  {45},  {London Mathematical Society Student Texts},   {Cambridge}.

\bibitem{CGC-cd-art-chen} Chen, Chin-Long, \textit{{Equivalent irreducible Goppa codes}},     IEEE Trans. Inf. Theory, 24, 766-770,  1978.

\bibitem{CGC-cd-art-patjohn1}  {Fitzpatrick, P. and Ryan, J.~A.},
 \textit{{Counting irreducible {G}oppa codes}},  {Journal of the Australian Mathematical Society},
 {2001}, {71},  {299--305}.

\bibitem{CGC-cd-art-patjohn2}
 {Fitzpatrick, P.  and Ryan, J.~A.},
 \textit{{The number of inequivalent irreducible {G}oppa codes}},
  {International Workshop on Coding and Cryptography, Paris, 2001}.






\bibitem{CGC-cd-phdthesis-marta}
 {Giorgetti, M.},
\textit{{On some algebraic interpretation of classical codes}},
 {University of Milan},
 {2006}.

\bibitem{CGC-alg-art-Guralnick1}
  {Guralnick, Robert M.},
  \textit{{Subgroups of prime power index in a simple group}},
  {Journal of Algebra},
  {81}, 2, 304--311
  {1983}.

\bibitem{MR1996953}
{Huffman, W. Cary and Pless, Vera},
\textit{{Fundamentals of error-correcting codes}},
{Cambridge University Press},
{Cambridge},
{2003}.


\bibitem{CGC-alg-book-hupper1}   {Huppert, B.}, \textit{{Endliche Gruppen  1}}, {Berlin, Heidelberg}, {Springer-Verlag},{ 1967}.



\bibitem{kantor_singer}
{Kantor, William M.},
\textit{{Linear groups containing a {S}inger cycle}},
{Journal of Algebra},
{62},
{1980},
{1},
{232--234}.






\bibitem{CGC-alg-Li}
 {Li, Cai Heng},
\textit{ {The finite primitive permutation groups containing an abelian
  regular subgroup}},
 {Proceedings of the London Mathematical Society. Third Series},
 {87},
 {2003}.

\bibitem{conjugacyaffine}
{Liebeck, Martin W. and Shalev, Aner},
\textit{{Maximal subgroups of symmetric groups}},
{Journal of Combinatorial Theory. Series A},
{75},
{1996},
{2},
{341--352}.


\bibitem{CGC-cd-art-McEli1} {McEliece, R.J.}, \textit{{A public key cryptosystem based on algebraic coding theory}}, {JPL DSN}, 114--116, {1978}.



 \bibitem{CGC-cd-book-macwilliamsI}
 {MacWilliams, F.~J. and Sloane, N.~J.~A.}, \textit{{The theory of error-correcting codes {I}}},
 {North-Holland Publishing Co.},
 {1977}.




\bibitem{CGC-cd-phdthesis-john}
 {Ryan, J.},
\textit{{Irreducible Goppa codes}}, Ph.D. Thesis,
 {University College Cork, Cork, Ireland},
 {2002}.
 \end{thebibliography}

\end{document}